\documentclass[12pt]{article}
\usepackage{epsfig}
\newcommand\C{{\bf C}}
\newcommand\e{{\bf e}}
\renewcommand\v{{\bf v}}
\newcommand\x{{\bf x}}
\newcommand\eref[1]{$(\ref{#1})$}
\renewcommand\Re{\mathop{\rm Re}}
\renewcommand\Im{\mathop{\rm Im}}
\newtheorem{conjecture}{Conjecture}
\title{A conjecture that the roots of a univariate
polynomial lie in a union of annuli  \\
(Interim Revised Version)\thanks{Supported in part by
NSF award 0434338.}}
\author{Stephen A.~Vavasis\thanks{Department of Computer Science,
4130 Upson Hall, Cornell University, Ithaca, NY 14853 USA, 
vavasis@cs.cornell.edu.}}
\begin{document}
\maketitle
\begin{abstract}
We conjecture that the roots of a degree-$n$ univariate
complex polynomial
are located in a union of $n-1$ annuli, each of which is
centered at a root of the
derivative and whose radii depend on higher
derivatives.  We prove the conjecture for the cases
of degrees 2 and 3,
and we report on tests with randomly generated polynomials of
higher degree.

We state two other closely related conjectures concerning Newton's
method.  If true, these conjectures imply the existence of a simple,
rapidly convergent algorithm for finding all
roots of a polynomial.
\end{abstract}

\section{Conjecture concerning annuli}
\label{sec:conj}
Let $p(z)$
be a univariate polynomial with coefficients in $\C$.  Let $z_1,\ldots,
z_n$ be its roots.  Let $\zeta_1,\ldots,\zeta_{n-1}$ be the roots
$p'$.  This paper proposes the conjecture
that  $z_1,\ldots,z_n$ lie in a union of $n-1$ annuli,
one for each of $\zeta_1,\ldots,\zeta_{n-1}$.  The two radii of each annulus
are determined from higher derivative values, and the inner radius
is a constant fraction of the outer radius.
The formal statement is as follows.

\begin{conjecture}
There exist two universal constants $0<\iota_1\le 1 \le \iota_2$
with the following property.
Let $z_1,\ldots, z_n$ be the roots of a degree-$n$ complex
univariate polynomial $p(z)$.  Let
$\zeta_1,\ldots,\zeta_{n-1}$ be the roots of its derivative.
Let $\rho_1,\ldots,\rho_{n-1}$ be defined by
\begin{equation}
\rho_j=\min_{k=2,\ldots,n} 
   \left| \frac{p(\zeta_j)k!}{p^{(k)}(\zeta_j)}\right|^{1/k}.
\label{eq:rhodef}
\end{equation}
Define annulus 
\begin{equation}
A_j=\{z: \iota_1\rho_j \le |z-\zeta_j| \le \iota_2\rho_j\}
\label{eq:Aidef}
\end{equation}
for $j=1,\ldots,n-1$.
Then for each $i=1,\ldots,n$,
\begin{equation}
z_i\in A_1 \cup \cdots\cup A_{n-1}.
\label{containment}
\end{equation}
\end{conjecture}

\begin{figure}
\begin{center}
\epsfig{file=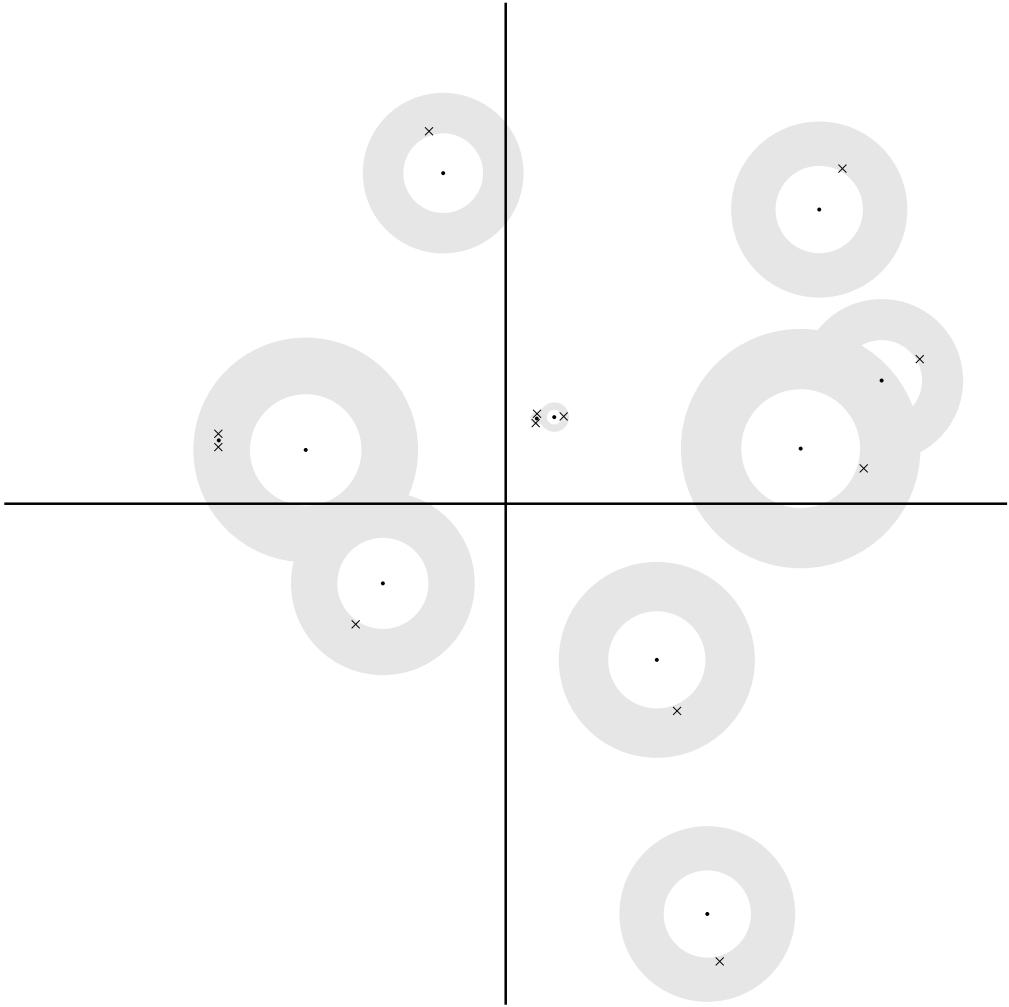, height=14cm}
\end{center}
\caption{The `$\times$' marks shown in the figure are the 12 roots of a
degree-12 polynomial plotted in the complex plane.  There is a cluster
of two very close roots and another cluster of three. 
The `.' markers show
the 11 roots of the derivative.  The $j$th derivative root
is surround
by an annulus with inner radius $.66\rho_j$ and outer radius
$1.33\rho_j$, where $\rho_j$ is given by \protect\eref{eq:rhodef}
and the scalars $.66$ and $1.33$ are chosen based on experiments in
Section~\protect\ref{sec:compu} as approximations to $\iota_1$ and $\iota_2$.}
\label{fig:example-poly}
\end{figure}

An example of the conjecture for a particular polynomial is
illustrated in Fig.~\ref{fig:example-poly}.

\begin{center}
$\cdots$
\end{center}

{\bf NOTE ADDED IN REVISION.}
Recently, M.~Giusti, J.~Heintz, 
G.~Lecerf, L.~Pardo,  B.~Salvy
and J.-C.~Yakoubsohn have shown
(unpublished communication) that
Conjecture 1 is partly true and
partly false.  In particular,
they have shown that $\iota_1$ 
exists and may be taken to be
$(\sqrt{5}-1)/2$.  They have also
constructed a family of polynomials
showing that $\iota_2$ cannot exist
(i.e., 
$|z-\zeta_j|/\rho_j\rightarrow\infty$)
as the degree tends to infinity for this family.
Their counterexample apparently does not
invalidate either Conjecture 2 or Conjecture 3
below.  A more detailed revision of this paper
is forthcoming.

\begin{center}
$\cdots$
\end{center}

This conjecture, if true, suggests a
simple algorithm based on Newton's method for finding all the roots
of a complex univariate polynomial.  
The algorithm (as well as two other related conjectures)
is described in Section~\ref{sec:newt}.

Observe that some of the terms in the min in \eref{eq:rhodef} could
be infinite, but they cannot all be infinite since $p^{(n)}$ is a nonzero
constant.
Observe also that the conjecture is invariant
under change of variables of the form $z'=az+b$ for any nonzero 
$a\in \C$ and any $b\in \C$.  It is also invariant under conjugation
and rescaling of the
polynomial (i.e., replacing $p$ by $cp$ or $c\bar{p}$ 
for a nonzero $c\in \C$).

The formula \eref{eq:rhodef} is motivated by considering the polynomial
$p(z)=z^n-c$.  The roots of this polynomial are the $n$ roots of $c$.
The derivative roots are all at 0.  For this polynomial, all
the terms in \eref{eq:rhodef} are infinite except the $k=n$ term,
and therefore $\rho_j=|c|^{1/n}$ for all $j$.  Thus, we could take
$\iota_1=\iota_2=1$ for this restricted class of polynomials.
Some additional background is given in Section~\ref{sec:newt}.
In Section~\ref{sec:deg3} we establish the conjecture for degree
2 and degree 3 polynomials.  In Section~\ref{sec:compu} we describe
a computational experiment to test the conjecture for high-degree
polynomials.

Aside from the application to Newton's method, the issue of location
of roots with respect to the derivative roots is an inherently interesting
matter that has attracted quite a bit of attention in the literature.
A classic result in this regard is the Gauss-Lucas theorem, which
states that the roots of $p'$ lie in the convex hull of the roots of
$p$.  The most comprehensive treatment of the relationships between
the roots of $p'$ and roots of $p$ appears to be the monograph of
\cite{sch}.  The present conjecture, however, appears to provide
more precise information about the location of the
roots of $p$ in terms of derivative roots than any of the
known theorems in \cite{sch}.

\section{Polynomials of degree 2 or 3}
\label{sec:deg3}
It follows immediately from the observations in the previous section
that the conjecture is true for polynomials of degree 2 since any
nondegenerate quadratic, after suitable change of variables and rescaling,
can be transformed to $z^2-1$.  Thus, for quadratics, the conjecture
is true with $\iota_1=\iota_2=1$.  
The conjecture is trivially true in the degenerate case of a
quadratic with a double root
since $\rho_1=0$ and $z_i-\zeta_1=0$.

In the case of degree-3 polynomials, first consider the cases
of a root with multiplicity 2 or 3.  For the multiplicity-3 case,
the conjecture is true for the same reason as in the last paragraph.
For the multiplicity-2 case, by rescaling and changing variables,
we may assume the polynomial is $p(z)=z^2(z-1)$.  In this case,
$z_1=z_2=0$ and $z_3=1$, while $\zeta_1=0$ and $\zeta_2=2/3$.  
Also, one checks that $\rho_1=0$ and $\rho_2=\min(\sqrt{4/27},(4/27)^{1/3})
\approx .3849$.  Since $|z_3-\zeta_2|=1/3$, for the conjecture to be
true in this case requires an $\iota_1\le (1/3)/(\sqrt{4/27})\approx
0.866$ and $\iota_2=1$.

The remainder of this section considers the case of a
cubic with three distinct roots.  Again by rescaling, changing
variables, and taking complex
conjugates, we may assume that the polynomial is monic and that
its closest pair of roots
lie at $-1$ and $1$, and that the
third root $a$ lies in Quadrant I.  Thus, $p(z)=(z-1)(z+1)(z-a)$ where
$\Re a\ge 0$, $\Im a\ge 0$, and $|a-1|\ge 2$. (The latter 
inequality follows by the assumption
that $-1$ and $1$ are closer to each other than to $a$.)
One checks
that $p'(z)=3z^2-2az-1$ hence
$$\zeta_1,\zeta_2=\frac{a\pm\sqrt{a^2+3}}{3}.$$
Since $a^2+3$ is in the upper half-plane by assumptions made about
$a$, we will assume that the branch of square-root appearing in the
preceding formula is chosen so that $\sqrt{a^2+3}$ lies in Quadrant I.
We will let $\zeta_1$ be the `$-$' branch and $\zeta_2$ the `$+$' branch.

Let $\zeta_\nu$ stand for either $\zeta_1$ or $\zeta_2$.  
Solving $p'(\zeta_\nu)=0$
yields the identity $a=(3\zeta_\nu^2-1)/(2\zeta_\nu)$ 
and hence $\zeta_\nu-a=
(-\zeta_\nu^2+1)/(2\zeta_\nu)$.

Let us now consider the three roots $\{z_1,z_2,z_3\}=\{1,-1,a\}$ in order.
For $z_1=1$, consider
the annulus about $\zeta_1$.
By definition, $\iota_1$ and $\iota_2$ must be
universal lower and upper bounds on the quantity
$$\frac{\min(|2p(\zeta_1)/p''(\zeta_1)|^{1/2},
|p(\zeta_1)|^{1/3})}{|\zeta_1-1|},$$
which, after using the facts that 
$p(\zeta_1)=(\zeta_1-1)(\zeta_1+1)(\zeta_1-a)$ and
$\zeta_1-a=(-\zeta_1^2+1)/(2\zeta_1)$,
simplifies to
\begin{equation}
\min\left(\left|\frac{(\zeta_1+1)^2}{3\zeta_1^2+1}\right|^{1/2},
\left|\frac{(\zeta_1+1)^2}{2\zeta_1(\zeta_1-1)}\right|^{1/3}\right).
\label{eq:z1quant}
\end{equation}
To establish an upper bound on \eref{eq:z1quant}, observe that
$(3\zeta_1^2+1)-2\zeta_1(\zeta_1-1)=(\zeta_1+1)^2$, hence
\begin{equation}
|\zeta_1+1|^2 \le |3\zeta_1^2+1| + |2\zeta_1(\zeta_1-1)|.
\label{eq:zeta1bd}
\end{equation}
Let $\alpha^*$ be chosen as the real root of
\begin{equation}
(\alpha^*)^3-\alpha^*-1=0
\label{eq:alphasdef}
\end{equation}
This $\alpha^*$ is close to
$1.3247$.  It follows from \eref{eq:alphasdef} that
$(1+1/\alpha^*)^{1/2}$ and $(1+\alpha^*)^{1/3}$ are both
equal to $\alpha^*$.   

Turning back to \eref{eq:zeta1bd}, take two subcases
depending on the relative sizes of the two terms on the
right-hand side of \eref{eq:zeta1bd}.
Subcase 1 is that $|3\zeta_1^2+1|\le \alpha^*|2\zeta_1(\zeta_1-1)|$.
In this case, the right-hand side of \eref{eq:zeta1bd} is
at most $(1+\alpha^*)\cdot |2\zeta_1(\zeta_1-1)|$, and hence
the second term of \eref{eq:z1quant} is at most $(1+\alpha^*)^{1/3}$,
which is equal to $\alpha^*$.

Subcase 2 is that $|3\zeta_1^2+1|\ge \alpha^*|2\zeta_1(\zeta_1-1)|$.
In this case, the right-hand side of \eref{eq:zeta1bd} is
at most $(1+1/\alpha^*)\cdot|3\zeta_1^2+1|$ and hence
the first term of \eref{eq:z1quant} is at most $(1+1/\alpha^*)^{1/2}$,
which is equal to $\alpha^*$.  This establishes the upper bound of
$\alpha^*$ on \eref{eq:z1quant}.

For the lower bound, we use the following cruder 
argument.
Let $Q_I$ denote the first quadrant.  Observe that the function of $a$
given by $\zeta_1(a)=(a-\sqrt{a^2+3})/3$ is analytic in the interior of $Q_I$
with a singularity on the boundary (at $\sqrt{3}i$), and therefore
its real part is harmonic.  This means that the minimum value of the real
part of $\zeta_1$ is attained either for $a$ 
on the boundary of $Q_I$ or in the limit for
infinitely large $a$.  Observe that $\zeta_1=a(1-\sqrt{1-3/a^2})/3$, and
for very large $|a|$, $\sqrt{1-3/a^2}\approx 1-1.5/a^2$, hence 
$\zeta_1\approx 1/(2a)$ which
tends to $0$ for large $|a|$.  Thus, the minimum real part of $\zeta_1$
is attained on the boundary rather than $\infty$. 
For $a$ on the positive imaginary axis (one boundary
of $Q_I$),
say $a=ti$, we compute that $\zeta_1=(-\sqrt{3-t^2}+ti)/3$, which has
real part equal to 0 if $t\ge \sqrt{3}$ else real part equal to 
$-\sqrt{3}/3$ or greater for $t\in[0,\sqrt{3}]$. Along the positive
real axis (the other boundary of $Q_I$), $\zeta_1$ is increasing, as one
can check from the derivative, hence the minimum value of the real
part is again at 0 and is equal to $-\sqrt{3}/3$.  Therefore,
$\Re(\zeta_1+1)\ge 1-\sqrt{3}/3$ hence 
\begin{equation}|\zeta_1+1|\ge 1-\sqrt{3}/3\approx 0.423.
\label{eq:zeta1abd}
\end{equation}
(This bound would be improved if we also accounted for the constraint
that $|a-1|\ge 2$.)

Next we have the following chain of inequalities to analyze the first term
of \eref{eq:z1quant}:
\begin{eqnarray}
|3\zeta_1^2+1| & = & |3\zeta_1^2+6\zeta_1+3 -6\zeta_1-6+4| \nonumber \\
  &\le & |3\zeta_1^2+6\zeta_1+3| + |6\zeta_1+6|+4 \nonumber \\
  &= & 3|\zeta_1+1|^2 + 6|\zeta_1+1| + 4 \nonumber\\
  &= &  3|\zeta_1+1|^2 + 6\frac{|\zeta_1+1|^2}{|\zeta_1+1|} + 4 \nonumber\\
  &\le & [3+6/(1-\sqrt{3}/3)+4/(1-\sqrt{3}/3)^2]\cdot |\zeta_1+1|^2,
\label{eq:3zbd}
\end{eqnarray}
where we have used the inequality \eref{eq:zeta1abd} to obtain
the last line.  The quantity in square brackets is approximately 39.6.

The second term can be similarly analyzed:
\begin{eqnarray}
|2\zeta_1(\zeta_1-1)|  & = & |2\zeta_1^2+4\zeta_1+2-6\zeta_1-6+4| \nonumber\\
 & \le & 2|\zeta_1+1|^2+6|\zeta_1+1|+4 \nonumber \\
& \le & [2+6/(1-\sqrt{3}/3)+4/(1-\sqrt{3}/3)^2]\cdot |\zeta+1|^2. 
\label{eq:2zbd}
\end{eqnarray}
Combining \eref{eq:3zbd} and \eref{eq:2zbd}
 establishes a rather poor lower bound of
$\sqrt{1/39.6}$ on \eref{eq:z1quant}.

For $z_2=-1$, we must obtain lower and upper bounds on 
$$\frac{\min(|2p(\zeta_1)/p''(\zeta_1)|^{1/2},
|p(\zeta_1)|^{1/3})}{|\zeta_1+1|},$$
which simplifies to
\begin{equation}
\min\left(\left|\frac{(\zeta_1-1)^2}{3\zeta_1^2+1}\right|^{1/2},
\left|\frac{(\zeta_1-1)^2}{2\zeta_1(\zeta_1+1)}\right|^{1/3}\right).
\label{eq:z2quant}
\end{equation}
To establish a upper bound on \eref{eq:z2quant}, use
an argument analogous to the preceding analysis of $z_1$:
$(3\zeta_1^2+1)-2\zeta_1(\zeta_1+1)=(\zeta_1-1)^2$, hence
$$|\zeta_1-1|^2 \le |3\zeta_1^2+1| + |2\zeta_1(\zeta_1+1)|.$$
Depending on the relative sizes of the
terms on the right-hand side,
either $|\zeta_1-1|^2\le (1+1/\alpha^*)\cdot |3\zeta_1^2+1|$, implying that 
the first term of \eref{eq:z2quant} is at most $\alpha^*$,
or else $|\zeta_1-1|^2\le (1+\alpha^*)\cdot |2\zeta_1(\zeta_1+1)|$,
implying that the second term of \eref{eq:z1quant} is at most $\alpha^*$.
Here, $\alpha^*$ was defined by \eref{eq:alphasdef}.

To establish a lower bound, again similar arguments are used.
We first claim that $\Re \zeta_1 \le 0$; this follows again
by considering the extremal cases for $\Re \zeta_1$ as a function
of $a$ as above.  This implies that $|\zeta_1-1|\ge 1$. 

Then we have the following chain of inequalities to analyze the first term
of \eref{eq:z2quant}:
\begin{eqnarray*}
|3\zeta_1^2+1| & = & |3\zeta_1^2-6\zeta_1+3 +6\zeta_1-6+4| \\
  &\le & |3\zeta_1^2-6\zeta_1+3| + |6\zeta_1-6|+4 \\
  &= & 3|\zeta_1-1|^2 + 6|\zeta_1-1| + 4 \\
  &= &  3|\zeta_1-1|^2 + 6\frac{|\zeta_1-1|^2}{|\zeta_1-1|} + 4 \\
  &\le & (3+6+4)|\zeta_1-1|^2,
\end{eqnarray*}
where we have used the inequality  $|\zeta_1-1|\ge 1$ to obtain
the last line. 

The second term can be similarly analyzed:
\begin{eqnarray*}
|2\zeta_1(\zeta_1+1)|  & = & |2\zeta_1^2-4\zeta_1+2+6\zeta_1-6+4| \\
 & \le & 2|\zeta_1-1|^2+6|\zeta_1-1|+4 \\
& \le & (2+6+4)|\zeta-1|^2.
\end{eqnarray*}

The last root to analyze is $z_3=a$; for this root we will consider
the annulus about derivative root $\zeta_2$ instead of $\zeta_1$.
The quantity to analyze is
$$\frac{\min(|2p(\zeta_2)/p''(\zeta_2)|^{1/2},
|p(\zeta_2)|^{1/3})}{|\zeta_2-a|},$$
which, after simplification, is equal to
\begin{equation}
\min\left(\left|\frac{4\zeta_2^2}{3\zeta_2^2+1}\right|^{1/2},
\left|\frac{4\zeta_2^2}{\zeta_2^2-1}\right|^{1/3}\right).
\label{eq:z3quant}
\end{equation}
An upper bound on \eref{eq:z3quant} is obtained by observing
that $4\zeta_2^2=3\zeta_2^2+1+\zeta_2^2-1$, hence
$$|4\zeta_2^2|\le |3\zeta_2^2+1|+|\zeta_2^2-1|.$$
Then, using the same logic as in the previous two cases, we conclude
that \eref{eq:z3quant} has $\alpha^*$ as an upper bound.

For the lower bound, observe that $\sqrt{a^3+3}$ lies in Quadrant
I (denoted $Q_I$) provided that $a\in Q_I$.  If $w_1,w_2$ are
any two complex numbers both lying in
$Q_I$, then $|w_1+w_2|\ge |w_1|$. Therefore, since 
$\zeta_2=(a+\sqrt{a^3+3})/3$, we conclude that $|\zeta_2|\ge |a|/3$.
The assumptions $|a-1|\ge 2$ and $a\in Q_I$ together imply that
$|a|\ge \sqrt{3}$, and therefore $|\zeta_2|\ge \sqrt{3}/3$, hence
$|\zeta_2^2|\ge 1/3$.

The inequality derived in the last paragraph yields a lower
bound on both terms of \eref{eq:z3quant}.  For the first term,
$|3\zeta_2^2+1|\le 3|\zeta_2^2|+1 \le 6|\zeta_2^2|.$  Therefore,
the first term of \eref{eq:z3quant} is at least $(2/3)^{1/2}$.
For the second term, a similar use of the previous
paragraph shows $|\zeta_2^2-1|\le 4|\zeta_2^2|$, and hence the
second term is at least 1.

This concludes the analysis of the $n=3$ case.  We have shown
an upper bound of $\alpha^*$ for $\iota_2$ and a lower bound of
$\sqrt{1/39.6}$ for $\iota_1$.  We have written a Matlab program
that computes $\iota_1$ and $\iota_2$ for each cubic polynomial 
with roots at $-1,1,a$,
where $a$ ranges over
a fairly dense
grid lying in the set $\{a\in Q_I: |a-1|\ge 2\}$.  We found
that $\iota_1$ appears to be approximately $0.82$ while $\iota_2$
appears to be exactly $\alpha^*$, and in particular, the upper bound
of $\alpha^*$ on \eref{eq:z3quant} appears to be tight.

\section{Computational experiment with higher degree polynomials}
\label{sec:compu}

In this section we describe our Matlab computational experiment  
with higher degree polynomials.

We experiment with three degrees: $n=10,20,40$.
For each $n$, we generate 3000 random polynomials.  Each polynomial
is chosen by selecting its $n$ roots uniformly at random on the unit
circle.  The rationale for this choice (as opposed to a distribution over
a 2-dimensional domain) is to greatly increase the likelihood of nearby
or clustered roots, which is a more difficult case for root-finding.

The roots of the derivative polynomial are then computed, as are
the parameters $\rho_1,\ldots,\rho_{n-1}$ given by \eref{eq:rhodef}.
For each polynomial root $z_i$, $i=1,\ldots,n$,
the program seeks the $j$ in $\{1,\ldots,n-1\}$ such that
$|z_i-\zeta_j|/\rho_j$ is closest to $1$.  Call this quotient $\iota_{z_i}$.
The program then tabulates the minimum and maximum $\iota_{z_i}$
encountered among all 3000 trials; these are taken to be  estimates
for $\iota_1$ and $\iota_2$.

The results are as follows.  For $n=10$, $\iota_1=0.67$ and
$\iota_2=1.32$; for $n=20$, $\iota_1=0.66$ and $\iota_2=1.33$;
for $n=40$, $\iota_1=0.66$ and $\iota_2 =1.33$.  Thus, there seems
to be little appreciable change in the experimental values of $\iota_1$
or $\iota_2$ as the degree increases.

We mention two subtleties concerning the implementation of this computational
experiment.
As mentioned above, the program selects the roots of the polynomial at
random on the unit circle and then computes the derivative roots.
The naive method to compute derivative roots, namely, form the standard
monomial basis for the polynomial, differentiate it term by term, and then
use the Matlab {\tt roots} function on the derivative, is unstable for
polynomials with clusters of roots.  This naive 
implementation gave incorrect experimental
results.  We found that a better method for finding derivative roots
is to compute the eigenvalues of the $(n-1)\times (n-1)$ matrix
$$M=\mbox{diag}(z_2,\ldots,z_n)-\e\v^T/n$$
where $\e$ is the vector of all 1's and $\v$ is the vector whose
$i$th entry is $z_{i+1}-z_1$.  A brief explanation of why these eigenvalues
are derivative
roots is as follows.  If $M\x=\zeta\x$, then $z_ix_i-d=\zeta x_i$ for
each $i=2,\ldots,n$, where $d=\v^T\x/n$.  Solving yields
$x_i=d/(z_i-\zeta)$.  Substituting this formula for $x_i$ into
$d=\v^T\x/n$ and simplifying yields
$1/(\zeta-z_1)+\cdots+1/(\zeta-z_n)=0$, which is
the same as $p'(\zeta)=0$.

It is easy to check that this method for derivative
roots works much better than the naive method
for contrived examples of polynomials with
root clusters, e.g., the polynomial $(z-2)^{10}-.01^{10}=0$, which has
10 roots lying on a circle of radius $0.01$ about the point 2.  We have
not, however, attempted a formal proof of stability of this method.

A second stability
subtlety is the computation
of $\rho_j$.  The naive method, namely, forming $p$ in
the standard monomial basis and differentiating term by term to obtain
all the derivatives in \eref{eq:rhodef}, is unstable.  Our implementation
uses the following method,
which appears to be stabler.  
Form the standard monomial representation of the 
polynomial $q(z)=p(z+\zeta_j)$ by multiplying together its
degree-1 factors $(z-z_1+\zeta_j)$, \ldots
$(z-z_n+\zeta_j)$ (or equivalently, by
applying the Matlab {\tt poly} function to the $n$ shifted roots 
$z_1-\zeta_j$, \ldots, $z_n-\zeta_j$).
Then derivatives of the form $p^{(k)}(\zeta_j)$ are directly obtained
from the coefficients of $q$.

\section{Application to Newton's method}
\label{sec:newt}

Recall that Newton's method for finding a root of a complex function
$\phi$ is given by the iteration
$$x^{k+1}=x^k-\phi(x^k)/\phi'(x^k).$$
Newton's method is known to converge quadratically to a nondegenerate
root $x^*$ if the starting point $x^0$ is sufficiently close to $x^*$.
The {\em basin of attraction} for $x^*$ is the set of starting
points $x^0$ such that
Newton's method will converge to $x^*$ for that starting point.  
Note that although Newton's method is asymptotically quadratic, a
point in the basin of attraction could lead to a sequence of iterates
that meanders far away from $x^*$ for
an arbitrary number of iterations before
asymptotic quadratic convergence takes hold.

This difficulty leads us to define the
{\em basin of fast convergence} for $x^*$ to be the set of $x^0$ such
that the sequence of iterates
generated by Newton's method starting from $x^0$
converges quadratically immediately (rather than asymptotically)
according to the following inequality:
$$|x^{k}-x^*| \le \left(\sqrt{1/2}\right)^{2^k-1}|x^0-x^*|$$
for all $k\ge 0$.  
This definition is similar to one from Blum et al.~\cite{Blum}.

For a degree-$n$ polynomial $p$ such that the roots of $p'$ are
$\zeta_1,\ldots,\zeta_{n-1}$, let us define the {\em DR-circles}
to be the set of circles about $\zeta_1,\ldots,\zeta_{n-1}$
of radius $\rho_1,\ldots,\rho_{n-1}$ respectively. (Here, DR stands
for ``derivative root.'')
By Conjecture 1, the roots of $p$ apparently all lie close to the union of
its DR-circles.  This suggests that many points on 
DR-circles lie in the basins of fast convergence of the roots
of $p$.  

There are several different possible conjectures that could be made
about DR-circles and basins of convergence.  The first conjecture is that
the basin of fast convergence
of each root of $p$ contains a subsegment of at least one DR circle.
Consideration of the polynomial $z^n-c$ indicates that the length
of this segment could be as small as $O(1/n)$ radians.
This is our conjecture:

\begin{conjecture}
There is a universal constant $\eta_1>0$ with the following
property.  Let $p$ be
a degree-$n$
univariate complex polynomial whose roots are $z_1,\ldots,z_n$.
Let $C_1,\ldots,C_{n-1}$ be the DR-circles of $p$, that is,
circles centered about the roots $\zeta_1,\ldots,\zeta_{n-1}$ of $p'$ such
that the radius of $C_j$ is $\rho_j$ as defined by \eref{eq:rhodef}.
Let $z_i$ be any root of $p$ of multiplicity $1$.  
Then there exists a segment of a DR-circle of length $\eta_1/n$
radians lying in the basin of fast convergence of $z_i$.
\end{conjecture}

To test this conjecture, 
we used the same set-up (3000 polynomials of degrees 10, 20 and 40)
as in Section~\ref{sec:compu}.  We discretized
each DR-circle with $10n$ evenly spaced points, and for each root
of $p$ and each DR-circle
we counted the number of such DR-circle points in the root's basin 
of fast convergence.  If the conjecture were true, this number would
always be greater than a positive constant $10\eta_1/(2\pi)$ for
at least one DR-circle.  In fact
the minimum number for $n=10$ was 7, for $n=20$ was 10, and for $n=40$
was 13.  This gives some evidence in favor of the conjecture.

If this conjecture were true, it would imply a very simple algorithm
based solely on Newton's method for finding all the roots of a 
degree-$n$ polynomial.
First, find the unique root of the linear polynomial $p^{(n-1)}$.  From this
root of $p^{(n-1)}$, find the two roots of $p^{(n-2)}$ by starting
Newton's method from a sufficient number of sample points
on the DR-circle of $p^{(n-2)}$.
Once the two roots of $p^{(n-2)}$ are found, construct the
DR-circles of $p^{(n-3)}$,
sample them with points, and carry out
Newton's method to find roots of $p^{(n-3)}$, etc., until finally
we find the roots of $p$ from those of $p'$. 
The conjecture suggests that the number
of sample points per DR-circle ought to be $O(n)$.

The complexity of this algorithm may be estimated as follows.
Suppose the roots lie in a disk of radius $R$ and root accuracy of
$\epsilon$ is desired.  Starting from the basin of fast convergence,
Newton's method requires $O(\log\log (R/\epsilon))$ iterations to
achieve the desired accuracy.  See Renegar \cite{Renegar} for a more
careful explanation of the factor $\log\log(R/\epsilon)$ as well as a
matching lower bound.
For finding the roots of $p^{(n-k)}$,
we require Newton's method to be started on $k-1$ circles, with
$O(k)$ points per circle.  Each iteration of Newton's method requires
$O(k)$ arithmetic operations.  Thus, the number of operations for the
roots of $p^{(n-k)}$ is $O(k^3\log\log(R/\epsilon))$.  
This is summed from $k=1,\ldots,n$,
yielding a bound of $O(n^4\log\log(R/\epsilon))$ operations.

The previous literature has many algorithms for finding all roots
of a univariate polynomial; see e.g., the survey
of Pan \cite{Pan}.  
Our complexity bound is worse than published bounds for rootfinding
in terms of its dependence on $n$,
although it is much simpler than most algorithms.  

Another rootfinding algorithm that uses only Newton's
method is due to Hubbard et al.~\cite{Hubbard} and is even simpler
than ours in that it uses a fixed set of Newton starting
points that depends only on the degree $n$.  
The drawback of the algorithm of Hubbard et al.\
is that Newton's method is in general 
not necessarily initiated in the basin
of fast convergence, so complexity estimates are far from optimal.

Renegar also has an algorithm \cite{Renegar} for all roots
of a univariate polynomial based primarily on Newton's method.
Renegar's algorithm always initiates Newton's
method in the basin of fast convergence and hence also has a running
time proportional to $\log\log(R/\epsilon)$ but has a better dependence
on $n$ than ours.  Renegar's algorithm, which
to some extent motivated the present work, is based on the following
key idea.  A root
$z$ of a polynomial $p$ has a large
basin of fast convergence unless $z$ is part of a cluster of closely
spaced roots.  Suppose, for
example, that three roots of $p$ are clustered, and no other root
is nearby.
In this case, all three will have
small basins of fast convergence.  On the other hand, it is 
guaranteed in this case that there is a root of $p''$ close to the three
clustered roots, and this root of $p''$ will have a large basin of fast
convergence.  More generally, an isolated cluster of $k$ roots must be
near a root of $p^{(k-1)}$  that has a large basin of fast convergence.
Thus, Renegar's algorithm consists of zooming in on 
root clusters (possibly recursively, since clusters can be nested
inside other clusters) by using Newton's method for roots of derivatives.
Once a point near the cluster is found, subdivision is used to find
the basin of fast convergence for each individual root in the cluster.
The drawback of Renegar's method is that, in addition to Newton's method,
it involves some other operations such as computation
of approximate winding numbers that might be difficult to implement
in practice.

A proof of Conjecture 2 might lead to further insight that would reduce
the $n^4$ factor in the complexity bound.
For example, suppose it were possible to predict 
which DR-circle would have at least one constant-sized segment
in the  union of basins of fast convergence.
In this case, we could modify
the above procedure by tracking only one DR-circle per derivative
and sampling
that circle with a constant number of points.  This
yields an $O(n^2\log\log(R/\epsilon))$ algorithm for finding a single
root of $p$.  
Then this root could be used to deflate the polynomial, and the
process could be repeated, yielding an $O(n^3\log\log(R/\epsilon))$
algorithm to find all roots.

We check this latter possibility by computing, for each polynomial
in our test set, what is the minimum number of DR-circles that have 
at least $1/10$ of their sample points (i.e., a total segment length of
$2\pi/10$ radians) in a basin of fast convergence.  This number appears
to grow linearly: for $n=10$, there were always at least 6 such DR-circles,
for $n=20$ there were always at least $13$, and for $n=40$ there were at
always at least $25$.  Let us
state this as another conjecture.

\begin{conjecture}
There exist two universal constants $\eta_2>0$, $\eta_3>0$ as follows.
For any degree-$n$ polynomial $p$, at least $\eta_2n$ of its DR-circles contain
segments of length $\eta_3$ radians that lie in the union of the
basins of fast convergence of roots of $p$.
\end{conjecture}

A final possible conjecture concerns the total length of basins
of convergence.  The total length of DR-circles in radians
is $2\pi(n-1)$, and $p$ has $n$ roots, so one might conjecture that 
the basin of fast convergence of a particular root meets DR-circles in at least
a constant number of total radians.  Our computational experiment, however,
did not support this conjecture---in fact, our experiments suggests that
it is more likely that the minimum number of total radians in the basin
of convergence for a particular root is $O(1/n)$, i.e., no better
than what Conjecture 2 implies for a single DR-circle.

Again, there is one subtle numerical stability issue concerning the
tests in this section.
If the Newton update term $p(z)/p'(z)$ is computed by writing down $p$
and $p'$
in standard monomial form and substituting the current iterate, then
an unstable procedure results and the computational
test yields false results.
Instead, the computational test
uses the mathematically equivalent formula
$$p(z)/p'(z)=\frac{1}{(z-z_1)^{-1}+\cdots+(z-z_n)^{-1}}.$$
This formula is applicable only in the case that the roots 
of $p$ are
already known, which obviously
would not occur in practical application of
Newton's method.

\section{Concluding remarks}
This paper raises three conjectures concerning the location of roots
of a polynomial.  Clearly, the main topic for future work would be
to prove one of them.  Another interesting question for practical
work concerns a numerically stable implementation of the root-finding
procedure outlined in the previous section.  We imagine that for numerical
stability, Newton's method should always be applied to shifted polynomials
(i.e., polynomials of the form $q(z)=p(z-s)$, where $p$ is the original
polynomial),
where the shift is close to the sought-after root.

A final interesting question is whether the conjecture generalizes to 
arbitrary entire functions whose roots have finite multiplicities.  
In this case,
$k$ appearing in \eref{eq:rhodef} ought to run from $2$ to $\infty$.

\section{Acknowledgment}
The author thanks James Renegar and Alex Vladimirsky of Cornell
for helpful discussions about this work.

\bibliography{univar}
\bibliographystyle{plain}

\end{document}